\documentclass{amsart}

\usepackage{amsmath}
\usepackage{amsthm}
\usepackage{amssymb,amsfonts}
\usepackage{euscript}

\usepackage[ps,matrix,arrow,curve,cmtip]{xy}

\title{Fibrations and homotopy colimits of simplicial sheaves}
\author{Charles Rezk}
\address{Department of Mathematics,
Northwestern University,
Evanston, IL 60208}
\email{ rezk@math.nwu.edu}
\date{November 3, 1998}

\numberwithin{equation}{section}
\makeatletter
  \let\c@subsection\c@equation
\makeatother

\theoremstyle{plain}   

\newtheorem{theorem}[subsection]{Theorem}
\newtheorem{propos}[subsection]{Proposition}
\newtheorem{corol}[subsection]{Corollary}
\newtheorem{lemma}[subsection]{Lemma}

\theoremstyle{remark}
\newtheorem{remark}[subsection]{Remark}    
\newtheorem{exam}[subsection]{Example}

\DeclareMathOperator{\colim}{colim}
\DeclareMathOperator*{\coliml}{colim}
\newcommand{\op}{{\operatorname{op}}}
\newcommand{\ob}{{\operatorname{ob}}}

\newcommand{\ra}{\rightarrow}
\newcommand{\lra}{\longrightarrow}
\newcommand{\xra}{\xrightarrow}

\newcommand{\cat}[1]{{\operatorname{\mathbf{#1}}}}
\newcommand{\overcat}[2]{{(#1\downarrow #2)}}

\newcommand{\Set}{{\operatorname{\EuScript{S}et}}}
\newcommand{\sSet}{{\operatorname{\EuScript{S}}}}

\DeclareMathOperator{\hocolim}{hocolim}

\newcommand{\realiz}[1]{\lvert#1\rvert}

\newcommand{\pullback}[1]{\underset{#1}{\times}}
\newcommand{\union}[1]{\underset{#1}{\cup}}

\newcommand{\N}{\mathbb{N}}

\newcommand{\eev}{\wedge}

\newcommand{\dfn}{\textbf}

\def\noloc{\;{:}\,}

\newcommand{\forcepar}{\mbox{}\par}

\newcommand{\bool}[1]{{\mathcal{#1}}}
\newcommand{\Sh}{{\operatorname{Sh}}}
\newcommand{\Psh}{{\operatorname{Psh}}}
\newcommand{\topos}[1]{{\mathcal{#1}}}
\DeclareMathOperator{\Ex}{Ex}
\newcommand{\propsub}[1]{{\bar{\mathcal{P}}#1}}
\newcommand{\powob}[1]{{\mathcal{P}#1}}

\begin{document}

\begin{abstract}
We show that homotopy pullbacks of sheaves of simplicial
sets over a Grothendieck topology distribute over  homotopy colimits;
this generalizes a result of 
Puppe about 
topological spaces.  In addition, we show that inverse image functors
between categories of simplicial sheaves preserve homotopy pullback
squares. 
The method we use introduces the notion of a
sharp map, which is analogous  to the notion of a
quasi-fibration of 
spaces, and seems to be of independent interest.
\end{abstract}

\maketitle

\section{Introduction}

Dold and Thom \cite{dold-thom;unendliche-symmetrische} introduced a class of
maps called quasi-fibrations.  A map $f\colon X\ra Y$ of topological
spaces is called a \emph{quasi-fibration} if for each point $y\in Y$
the fiber $f^{-1}(y)$ is naturally weakly equivalent to the homotopy
fiber of $f$ over $y$.  Thus, quasi-fibrations behave for some purposes
of homotopy theory very much like other types of fibrations; for
example, there is a long exact sequence relating the homotopy groups of $X$,
$Y$, and $f^{-1}(y)$.  A notable feature of quasi-fibrations is that
(as shown by Dold and Thom) quasi-fibrations defined over the elements
of an open cover of a space $Y$ can sometimes be ``patched''
together to give a quasi-fibration mapping to all of $Y$.

In this paper we study a class of maps called sharp maps.  In our
context, a map
$f\colon X\ra Y$ will be called \emph{sharp} if for each base-change of
$f$ along 
any map into the base $Y$ the resulting pullback square is homotopy
cartesian.  

We are particularly interested in sharp maps of sheaves of simplicial sets.
We shall show that sharp maps of sheaves of simplicial sets have
properties analogous to those of quasi-fibrations of topological 
spaces.  In particular, they can be ``patched together'', in a sense
analogous to the way that quasi-fibrations can be patched together.
We give several applications.

\subsection{Applications}

Let $\topos{E}$ denote a Grothendieck topos; that is, a category
equivalent to a category of sheaves on a small Grothendieck site.  The
category $s\topos{E}$ of \emph{simplicial} objects in $\topos{E}$
admits a Quillen 
closed model category structure, as was shown by Joyal (unpublished),
and by Jardine in 
\cite{jardine;simplicial-presheaves;pure} and
\cite{jardine-boolean-localization}.

Let $X\colon\cat{I}\ra s\topos{E}$ be a diagram of simplicial sheaves
indexed on a small category $\cat{I}$.  We say that such a diagram is
a \dfn{homotopy colimit diagram} if the natural map
$\hocolim_{\cat{I}}X\ra\colim_{\cat{I}}X$ is a weak equivalence of
objects in $s\topos{E}$, where $\hocolim$ denotes the homotopy colimit
functor for simplicial sheaves, generalizing that defined by Bousfield
and Kan \cite{yellow-monster} for simplicial sets. 

Given a map $f\colon X\ra Y$ of $\cat{I}$-diagrams of simplicial
sheaves, for each object $i$ of $\cat{I}$ there exists a commutative square
\begin{equation}\label{eq-square-of-maps-to-colims}
\vcenter{\xymatrix{
{Xi} \ar[r] \ar[d]_{fi}
& {\colim_{\cat{I}}X} \ar[d]
\\
{Yi} \ar[r]
& {\colim_{\cat{I}}Y}
}}\end{equation}
and for each morphism $\alpha\colon i\ra j$ of
$\cat{I}$ there exists a commutative square
\begin{equation}\label{eq-square-of-maps-in-diagrams}
\vcenter{\xymatrix{
{Xi} \ar[r]^{X\alpha} \ar[d]_{fi}
& {Xj} \ar[d]^{fj}
\\
{Yi} \ar[r]^{Y\alpha} 
& {Yj}
}}\end{equation}

The following theorem essentially says that in a category of
simplicial sheaves, homotopy pullbacks ``distribute'' over homotopy
colimits. 
\begin{theorem}
\label{thm-hocolims}
Let $f\colon X\ra Y$ be a map of $\cat{I}$-diagrams of simplicial
objects in a topos $\topos{E}$, and suppose that $Y$ is a homotopy
colimit diagram.  Then the 
following two properties hold.
\begin{enumerate}
\item [(1)]
If each square of the form (\ref{eq-square-of-maps-to-colims}) is
homotopy cartesian, then $X$ is a homotopy colimit diagram.
\item [(2)]
If $X$ is a homotopy colimit diagram, and each diagram of the
form (\ref{eq-square-of-maps-in-diagrams}) is homotopy cartesian, then
each diagram of the form (\ref{eq-square-of-maps-to-colims}) is also
homotopy cartesian.
\end{enumerate}
\end{theorem}

The proof of (\ref{thm-hocolims}) is given in Section~\ref{sec-ho-colims}.
This result is well-known when $s\topos{E}$ is the category of
simplicial sets:  Puppe \cite{puppe-remark-homotopy-fibrations}
formulates and proves
a version of the above result for the category of topological
spaces, which can be used to derive (\ref{thm-hocolims}) for
simplicial sets; see \cite[Appendix HL]{farjoun-cellular-spaces} for
more discussion of Puppe's result.
Also, Chach\'olski~\cite{chacholski-closed-classes}
has proved a result of this type in the category of
simplicial sets using purely simplicial methods. 

As another application we give the following.  Let
$p\colon\topos{E}\ra\topos{E}'$ be a geometric morphism of
Grothendieck topoi, and let $p^*\colon \topos{E}'\ra \topos{E}$
denote the corresponding inverse image functor.  This functor prolongs to a
simplicial functor $p^*\colon s\topos{E}'\ra s\topos{E}$.

\begin{theorem}\label{thm-inv-image-ho-cartesian}
The inverse image functor $p^*\colon s\topos{E}'\ra s\topos{E}$
preserves homotopy cartesian squares. 
\end{theorem}
The proof of (\ref{thm-inv-image-ho-cartesian}) is given in
Section~\ref{sec-basic-props-of-sharp-maps}.  An example of an inverse
image functor is the sheafification functor $L^2\colon
\Psh{\cat{C}}\ra \Sh{\cat{C}}$ associated to a Grothendieck topology
on $\cat{C}$.  Thus,
(\ref{thm-inv-image-ho-cartesian}) shows in particular  that
sheafification functors 
preserve homotopy cartesian squares.

\subsection{Organization of the paper}

In Section~\ref{sec-sharp-maps} we define sharp maps and state some of
their general properties.  In Section~\ref{sec-topos-facts} we recall
facts about sheaf theory and the model category structure on
simplicial sheaves.  Section~\ref{sec-local-charac-of-sharp-maps}
gives several useful characterizations of sharp maps of simplicial
sheaves, which are used to prove a number of properties in
Section~\ref{sec-basic-props-of-sharp-maps}, as well as the proof of
(\ref{thm-inv-image-ho-cartesian}).

In Section~\ref{sec-diagonal-of-s-obj} we prove a result about how
sharp maps are preserved by taking the diagonal of a simplicial
object.  This result is used in Section~\ref{sec-ho-colims} to prove a
similar fact about how sharp maps are preserved by homotopy colimits;
this result is used in turn to give a proof of
(\ref{thm-hocolims}).  Section~\ref{sec-special-diagrams}
does the hard work of showing that sharp maps which agree ``up to
homotopy'' can be glued together, thus providing lemmas which were needed for
Section~\ref{sec-diagonal-of-s-obj}.  

Section~\ref{sec-sharp-maps-in-bool-loc} proves a result about sharp
maps in a boolean localization which was needed in
Section~\ref{sec-local-charac-of-sharp-maps}.

In Section~\ref{sec-local-fibs-are-global-fibs} we prove that in a
boolean localization the local fibrations are the same as the global
fibrations, a fact which is used at several places in this paper.

\subsection{Acknowledgements}

I would like to thank Phil Hirschhorn, Mike Hopkins, Mark Johnson, 
Brooke Shipley, 
and Carlos Simpson for their comments on various parts of this work.

\section{Sharp maps}
\label{sec-sharp-maps}


In this section, we define the notion of a sharp map in a general
closed model category, and prove some of its general properties.  I
learned about the notion of a ``sharp'' map from Mike Hopkins, who was
originally led, for different reasons, to formulate the dual notion of a
``flat'' map.

Let $\cat{M}$ be a closed model category \cite{homotopical-algebra},
\cite{quillen-ratl-homotopy}. 
We say that a map $f\colon X\ra Y$ in 
$\cat{M}$ is \dfn{sharp} if for each diagram in $\cat{M}$
of the form 
$$\xymatrix{
  {A} \ar[r]^i \ar[d]
  & {A'} \ar[r] \ar[d]
  & {X} \ar[d]^f
  \\
  {B} \ar[r]^j
  & {B'} \ar[r]
  & {Y}
}$$
in which $j$ is a weak
equivalence and each square is a pullback square, the map $i$ is also a
weak equivalence. 
It follows immediately from the definition that the class of sharp
maps is closed under base-change.

\subsection{Proper model categories}

A model category $\cat{M}$ is said to be \dfn{right proper} if for
each pullback diagram in $\cat{M}$ of the form
$$\xymatrix{
{X'} \ar[r]^i \ar[d]
& {X} \ar[d]^f
\\
{Y'} \ar[r]^j
& {Y}
}$$
such that $f$ is a fibration and $j$ is a weak equivalence, then $i$
is also a weak equivalence.  The categories of topological spaces and
simplicial sets are two well-known examples of right-proper model
categories.  

There is an dual notion, in which a model category for
which pushouts of weak equivalences along cofibrations are weak
equivalences is called \dfn{left proper}.  A model category is
\dfn{proper} if it is both left and right proper.

Since the class of fibrations in a model
category is closed under base-change, we have the following.
\begin{propos}
\label{prop-fibs-are-sharp}
A model category $\cat{M}$ is right proper if and only if each
fibration is sharp.
\end{propos}

\subsection{Homotopy cartesian squares}
\label{subsec-ho-cartesian-squares}

Let 
\begin{equation}\label{eq-pullback-square}
\vcenter{\xymatrix{ 
{P} \ar[r] \ar[d] 
& {Y} \ar[d]^g 
\\ 
{X} \ar[r]^f 
& {B} 
}}
\end{equation}
be a commutative square in $\cat{M}$.  Say such a square is \dfn{homotopy
cartesian} if for some choice of factorizations $X\ra X' \ra B$ and
$Y\ra Y'\ra B$ of $f$ and $g$ into weak equivalences followed by
fibrations, the 
natural map $P\ra X'\times_{B}Y'$ is a weak equivalence.  It is
straightforward to show that the choice of factorizations does not matter.

Clearly, any \emph{pullback} square of the form
(\ref{eq-pullback-square}) in which $f$ and $g$ 
are already fibrations is homotopy cartesian.  Any square
weakly equivalent to a homotopy cartesian square is itself homotopy
cartesian. 

\begin{lemma}
\label{lemma-one-fib-gives-ho-cartesian}
In a right proper model category, a pullback square as in
(\ref{eq-pullback-square}) in which $g$ is a fibration is a homotopy
cartesian square.
\end{lemma}
\begin{proof}
Choose a factorization $pi$ of $f$ into a weak
equivalence $i\colon X\ra X'$
followed by a fibration $p\colon X'\ra B$.  Then we obtain pullback
squares 
\begin{equation}\label{eq-factored-pb-square}
\vcenter{\xymatrix{
{P} \ar[r]^j \ar[d]
& {P'} \ar[r] \ar[d]^h
& {Y} \ar[d]^g
\\
{X} \ar[r]^i 
& {X'} \ar[r]^p
& {B}
}}
\end{equation}
in which $j$ is a weak equivalence by (\ref{prop-fibs-are-sharp}),
since it is obtained by pulling 
back the weak equivalence $i$ along the fibration $h$.  Thus the
square (\ref{eq-pullback-square}) is weakly 
equivalent to the right-hand square of (\ref{eq-factored-pb-square})
which is homotopy cartesian.
\end{proof}

The following proposition gives the characterization of sharp maps
which was alluded to in the introduction; it holds only in a right
proper model category.
\begin{propos}
\label{prop-all-pbs-ho-cart-is-sharp}
In a right proper model category, a map $g\colon Y\ra B$ is a sharp
map if and only if each pullback square
(\ref{eq-pullback-square}) is a homotopy
cartesian square.
\end{propos}
\begin{proof}
First suppose that $g$ is sharp.
As in the proof of (\ref{lemma-one-fib-gives-ho-cartesian})
choose a factorization $pi\colon X\ra X'\ra B$ of $f$ into a weak
equivalence $i$ followed by a fibration $p$, obtaining a diagram
(\ref{eq-factored-pb-square}).  Then the right hand square of this
diagram is homotopy cartesian by
(\ref{lemma-one-fib-gives-ho-cartesian}), and $i$ and $j$ are
weak equivalences, since $j$ is the pullback of the weak equivalence
$i$ along the sharp map $g$.

Conversely, suppose $g$ is a map such that each pullback along $g$ is
a homotopy cartesian square.  Given a diagram of pullback squares as
in (\ref{eq-factored-pb-square}) in which $i$ is a weak equivalence,
it follows that $j$ is also a weak equivalence, since both the
right-hand square and the outer rectangle are homotopy cartesian
squares which are weakly equivalent at the three non-pullback
corners.  Thus $g$ is sharp.
\end{proof}

\begin{exam}
The category of topological spaces is a right proper model category.
The class of sharp maps of topological spaces includes all Serre
fibrations, as well as all fiber bundles.  Every sharp map is clearly
a quasi-fibration in the sense of Dold and Thom
\cite{dold-thom;unendliche-symmetrische}.  
It is not the case that all quasi-fibrations are sharp; indeed, the
class of quasi-fibrations is not closed under base change, see
\cite[Bemerkung 2.3]{dold-thom;unendliche-symmetrische}.
I do not know of a simple characterization of sharp maps of
topological spaces.
\end{exam}

\section{Facts about topoi}
\label{sec-topos-facts}

In this section we recall facts about sheaves and simplicial sheaves.
Our main reference for sheaf theory is Mac Lane-Moerdijk
\cite{moerdijk-maclane;sheaves-geometry-logic}.  

\subsection{Grothendieck topoi}

A \dfn{Grothendieck topos} $\topos{E}$ is a category equivalent to
some category $\Sh\cat{C}$ of sheaves of sets on a small Grothendieck
site 
$\cat{C}$.  Among the many properties of a Grothendieck topos
$\topos{E}$, we note that $\topos{E}$ has all small limits and
colimits, and that $\topos{E}$ is cartesian closed.  The internal hom
object in $\topos{E}$ is denoted by $Y^X$.

\begin{exam}\forcepar
\begin{enumerate}
\item
The category $\Set$ is a Grothendieck topos, since it is sheaves on a
one-point space.
\item
The presheaf category $\Psh\cat{C}$, defined to be the category of
functors $\cat{C}^\op\ra\Set$, is the category of sheaves of sets in the
trivial topology on $\cat{C}$, and thus is a Grothendieck topos.
\item
The category $\Sh(T)$ of sheaves of sets on a topological space $T$ is
a Grothendieck topos.
\end{enumerate}
\end{exam}

A \dfn{geometric morphism} $f\colon \topos{E}\ra\topos{E}'$ is a pair
of adjoint functors
$$f^\ast\colon \topos{E}'\rightleftarrows\topos{E}\noloc f_\ast$$
such that the left adjoint $f^\ast$ preserves finite limits.  The left
adjoint $f^\ast$ is called the \dfn{inverse image} functor, and
$f_\ast$ the \dfn{direct image} functor.

\subsection{Boolean localizations}\label{subsec-boolean-localization}

Let $\bool{B}$ be a complete Boolean algebra.  Then $\bool{B}$, viewed
as a category via the partial order on $\bool{B}$, has a natural
Grothendieck topology, and hence gives us a Grothendieck topos
$\Sh\bool{B}$.  (This topos is discussed in more detail in
Section~\ref{sec-local-fibs-are-global-fibs}.) 

A \dfn{Boolean localization} of a topos $\topos{E}$ is a geometric
morphism $p\colon\Sh\bool{B}\ra\topos{E}$ such that the inverse image
functor $p^\ast\colon\topos{E}\ra\Sh\bool{B}$ is faithful.

\begin{exam}\label{exam-boolean-loc}\forcepar
\begin{enumerate}
\item The category of sets is its own Boolean localization, since it
is equivalent to sheaves on the trivial Boolean algebra.

\item
For a category $\cat{C}$, let $\cat{C_0}\subset\cat{C}$ denote the
subcategory consisting of all objects and all identity maps.  Then
$p\colon \Psh\cat{C_0}\ra\Psh\cat{C}$ is a boolean localization, where
$p^\ast\colon\Psh\cat{C}\ra\Psh\cat{C_0}$ is the obvious restriction
functor; this
is because $\Psh\cat{C_0}$ is equivalent to the category of sheaves
on the boolean algebra $\mathcal{P}(\ob\cat{C})$, the power set of
$\ob\cat{C}$. 

\item
For a topological space $T$, let $T^\delta$ denote the underlying  set
$T$ with 
the discrete topology.  Then
$\Sh(T^\delta)\approx\Sh(\mathcal{P}T^\delta)$ is a boolean
localization of $\Sh(T)$; the inverse image functor $p^\ast\colon
\Sh(T)\ra\Sh(T^\delta)$ sends a sheaf $X$ to the collection of all
stalks of $X$ over every point of $T$.
\end{enumerate}
\end{exam}

\begin{remark}
In each of the examples above, the Boolean localization turned out to
be equivalent to a product of copies of $\Set$.  However, there exist
topoi $\topos{E}$ which do not admit a Boolean localization of this
type.
\end{remark}

Boolean localizations have the following properties.
\begin{enumerate}
\item Every Grothendieck topos has a Boolean localization.

\item The inverse image functor $p^\ast$ associated to a Boolean
localization functor $p\colon \Sh\bool{B}\ra\topos{E}$ reflects isomorphisms,
monomorphisms, epimorphisms, colimits, and finite limits.

\item The topos
$\Sh\bool{B}$ has a ``choice'' axiom: every epimorphism in
$\Sh\bool{B}$ admits a section.

\item The topos $\Sh\bool{B}$ is \dfn{boolean}; that is, each
subobject $A\subset X$ in $\Sh\bool{B}$ admits a ``complement'', namely
a subobject $B\subset X$ such that $A\cup B= X$ and $A\cap
B=\varnothing$.
\end{enumerate}
Property (1) is shown in
\cite[IX.9]{moerdijk-maclane;sheaves-geometry-logic}.   
See Jardine~\cite{jardine-boolean-localization} for proofs of the
other properties.

\subsection{A distributive law}

For our purposes it is important to note the following relationship
between colimits and pullbacks in a topos $\topos{E}$.  

\begin{propos}\label{prop-distributive-law}
Let $Y\colon\cat{I}\ra\topos{E}$ be a functor from a small category
$\cat{I}$ to a topos, and let $A\ra B=\colim_{\cat{I}}(i\mapsto Yi)$ be
a map.  Then the natural map
$$\colim_{\cat{I}}(i\mapsto A\pullback{B}Yi)\ra A$$
is an isomorphism.
\end{propos}
This proposition says that if an object is pulled back along a colimit
diagram, then that object can be recovered as the colimit of the
pulled-back diagram.
It makes sense to think of this as a ``distributive law''.  In fact,
in the special case in which $B=X_1\amalg X_2$, and $A\ra B$ is the
projection $(X_1\amalg X_2)\times Y\ra X_1\amalg X_2$ the
proposition reduces to the usual distributive law of products over
coproducts: $(X_1\times Y)\amalg (X_2\times Y)\approx(X_1\amalg
X_2)\times Y$.

To prove (\ref{prop-distributive-law}), note that it is true if
$\topos{E}=\Set$, and thus is true if $\topos{E}=\Psh\cat{C}$.  The
general result now follows from the properties of the sheafification
functor $L^2\colon\Psh\cat{C}\ra\Sh\cat{C}$. 

\begin{remark}
Consider the diagram $X$ defined by $i\mapsto A\times_{B}Yi$.
It is equipped with a natural transformation $f\colon X\ra Y$ with the
property that for each $\alpha\colon i\ra j$ in $\cat{I}$ the map
$X\alpha$ is the pullback of $Y\alpha$ along $fj$.

One can formulate the following ``converse'' to
(\ref{prop-distributive-law}) which is false.  Namely, given
a natural transformation $f\colon X'\ra Y$ of $\cat{I}$-diagrams such
that for each $\alpha\colon i\ra j$ in $\cat{I}$ the map $X'\alpha$ is
the base-change of $Y\alpha$ along $fj$, one may ask
whether the natural maps $X'i\ra A\times_{B}Yi$ are isomorphisms,
where $A=\colim_{\cat{I}}X'$.
A counterexample in $\topos{E}=\Set$ is to take $\cat{I}$ to be
a group $G$ and $X'\ra Y$ to be any map of non-isomorphic $G$-orbits.

(\ref{thm-hocolims}, part 1) may be viewed as a homotopy theoretic
analogue of (\ref{prop-distributive-law}).
(\ref{thm-hocolims}, part 2) may be viewed as a homotopy theoretic
analogue of the ``converse'' to (\ref{prop-distributive-law}).
\end{remark}

\subsection{Simplicial sheaves}

We let $s\topos{E}$ denote the category of simplicial objects in a
Grothendieck topos $\topos{E}$.  Note that $s\topos{E}$ is itself a
Grothendieck topos.  The full subcategory of \dfn{discrete}
simplicial objects in $s\topos{E}$ is equivalent to $\topos{E}$; thus,
we regard $\topos{E}$ as a subcategory of $s\topos{E}$.

For any topos $\topos{E}$ there is a natural functor
$\Set\ra\topos{E}$ sending a set $X$ to the corresponding
\dfn{constant} sheaf. This prolongs to a functor $s\Set\ra
s\topos{E}$, and we will thus regard any simplicial set as a constant
simplicial sheaf.

\subsection{Model category for simplicial sheaves}
\label{subsec-model-cat-simplicial-sheaves}

We will make use of the elegant model category structure on simplicial
sheaves provided by Jardine in \cite{jardine-boolean-localization}.
We summarize here the main properties of this structure which we need.
Let
$s\topos{E}$ denote the category of simplicial objects in a topos
$\topos{E}$.  Let $p\colon\Sh\bool{B}\ra\topos{E}$ denote a fixed
boolean localization of $\topos{E}$.  A map $f\colon X\ra Y$ in
$s\topos{E}$ is said to be
\begin{enumerate}
\item
a \dfn{local weak
equivalence} (or simply, a \dfn{weak equivalence}) if 
$$(L^2\Ex^\infty p^\ast
f)(b)\colon (L^2\Ex^\infty p^\ast X)(b)\ra (L^2\Ex^\infty p^\ast Y)(b)$$
is a weak equivalence for each 
$b\in\bool{B}$.  Here $\Ex^\infty\colon s\Sh\bool{B}\ra s\Psh\bool{B}$
denotes the 
functor obtained by applying Kan's $\Ex^\infty$ functor
\cite{kan;c-s-s-complexes} at each $b\in\bool{B}$, and $L^2$ denotes
the simplicial prolongation $s\Psh\bool{B}\ra s\Sh\bool{B}$ of the
sheafification functor.
\item
a \dfn{local
fibration} if $p^\ast f(b)\colon p^\ast X(b)\ra p^\ast Y(b)$ is a Kan
fibration 
for each $b\in\bool{B}$.  It should be pointed out that local
fibrations are not in general the fibrations in the model category
structure on $s\topos{E}$; but note
(\ref{prop-local-fibs-are-global}).
\item
a \dfn{cofibration} if it is a monomorphism.
\item
a \dfn{global fibration} (or simply, a \dfn{fibration}) if it has the
right lifting property with respect to all maps which are both
cofibrations and weak equivalences.  
\end{enumerate}

Note that the definition of local weak equivalence simplifies when
$\topos{E}=\Sh\bool{B}$, since $\Sh\bool{B}$ is its own boolean
localization.  Furthermore, a map $f$ in $s\topos{E}$ is a local weak
equivalence if and only if $p^*f$ in $s\Sh\bool{B}$ is a local weak
equivalence.  

\begin{theorem}
(Jardine \cite{jardine-boolean-localization},
\cite{jardine;simplicial-presheaves;pure})
The category $s\topos{E}$ with the above classes of cofibrations,
global fibrations, and
local weak equivalences is a proper simplicial  closed model category.
Furthermore, the characterizations of local weak equivalences,
local fibrations, and global fibrations do not depend on the choice of
boolean localization. 
\end{theorem}

\begin{exam}\forcepar
\begin{enumerate}
\item  When $s\topos{E}=s\Set$ this model category structure coincides
with the usual one, and local fibrations coincide with global
fibrations.
\item  For $s\topos{E}=s\Psh\cat{C}$, a map $f\colon X\ra Y$ is a
local weak equivalence, cofibration, or local fibration if for each
$C\in\ob C$, the map $fC\colon
X(C)\ra Y(C)$ is respectively a weak equivalence, monomorphism, or Kan
fibration of simplicial sets.
\item  For $s\topos{E}=s\Sh(T)$, a map $f\colon X\ra Y$ is a local
weak equivalence, cofibration, or local fibration if for each point $p\in T$
the map $f_p\colon X_p\ra Y_p$ of stalks is respectively a weak
equivalence, monomorphism, or Kan fibration of simplicial sets.
\end{enumerate}
\end{exam}

We also need the following property.

\begin{propos}\cite[Lemma 13(3)]{jardine-boolean-localization}
\label{prop-model-cat-colim-lemma}
Let $\topos{E}$ be sheaves on a Grothendieck topos.  
If $f_i\colon X_i\ra Y_i$ is a family of local weak equivalences in
$s\topos{E}$ indexed by a set $I$, then the induced map $f\colon
\coprod_{i\in I}X_i\ra\coprod_{i\in I}Y_i$ is a local weak equivalence.
\end{propos}

We need one additional fact about fibrations in a boolean localization.
\begin{propos}
\label{prop-local-fibs-are-global}
In the category $s\Sh\bool{B}$ of simplicial sheaves on a complete
boolean algebra $\bool{B}$,
the local fibrations are precisely the global fibrations.
\end{propos}
The proof of (\ref{prop-local-fibs-are-global}) is given in
Section~\ref{sec-local-fibs-are-global-fibs}. 

Finally, we note that if $f\colon\topos{E}\ra \topos{E'}$ is a
geometric morphism, then the induced inverse image functor $f^*\colon
s\topos{E'}\ra s\topos{E}$ preserves cofibrations and weak
equivalences; this is because the composite 
$$s\topos{E'}\xra{f^*} s\topos{E}\xra{p^*} s\Sh\bool{B}$$
must preserve such if $p\colon\Sh\bool{B}\ra\topos{E}$ is a boolean
localization of $\topos{E}$.

\subsection{Model category for simplicial presheaves}
\label{subsec-presheaf-model-cat}

Although we will not make much use of it here, we note that if
$\topos{E}=\Sh\cat{C}$ for some Grothendieck site $\cat{C}$, then
Jardine \cite[Thm. 17]{jardine-boolean-localization} constructs a
``presheaf'' closed model
category structure on $s\Psh\cat{C}$ related to that on $s\topos{E}$
(and not to be confused with the ``sheaf'' model category structure
obtained by 
applying the remarks of the previous section to
$\topos{E}=\Psh\cat{C}$).   In this structure on $s\Psh\cat{C}$, the
cofibrations are the 
monomorphisms, and the weak equivalences are the maps in $s\Psh\cat{C}$
which sheafify
to local weak equivalences in $s\topos{E}$.  Furthermore, a map in
$s\Psh\cat{C}$ is called a local fibration if it sheafifies to a local
fibration in $s\topos{E}$.   The natural adjoint pair
$s\Psh\cat{C}\leftrightarrows s\topos{E}$ induces an equivalence of
closed model categories in the sense of Quillen; in particular, the
homotopy category of $s\Psh\cat{C}$ (induced by the presheaf model category
structure) is equivalent to the homotopy category of $s\topos{E}$.
Thus, many results stated for $s\topos{E}$ such as
(\ref{thm-hocolims}) carry over to the presheaf model category of
$s\Psh\cat{C}$ without change.

\section{Local character of sharp maps of simplicial sheaves}
\label{sec-local-charac-of-sharp-maps}

The following theorem provides several equivalent characterizations of
sharp maps in $s\topos{E}$.  There are two types of such statements:
(5) and (6) say that sharpness is a ``local condition'', i.e.,
sharpness is detected on boolean localizations, while (2), (3),
and (4) say that sharpness is detected on ``fibers'', i.e., by pulling
back to the product of a discrete object and a simplex.

\begin{theorem}
\label{thm-equiv-charac-of-sharp}
Let $f\colon X\ra Y$ be a map of simplicial objects in a Grothendieck
topos $\topos{E}$.  The following are equivalent.
\begin{enumerate}
\item [(1)] $f$ is sharp.

\item [(2)] For each $n\geq 0$ and each map $S\ra Y_n$ in
$\topos{E}$, the induced pullback square
$$\xymatrix{
  {P} \ar[r] \ar[d]
  & {X} \ar[d]^f 
  \\
  {S\times\Delta[n]} \ar[r] 
  & {Y}
}$$
is homotopy cartesian.

\item [(3)] For each $n\geq 0$ there exists an epimorphism $S_n\ra
Y_n$ in $\topos{E}$ such that the induced pullback square
$$\xymatrix{
  {P} \ar[r] \ar[d]
  & {X} \ar[d]^f 
  \\
  {S_n\times\Delta[n]} \ar[r] 
  & {Y}
}$$
is homotopy cartesian.

\item [(4)] For each $n\geq 0$ there exists an epimorphism $S_n\ra
Y_n$ in $\topos{E}$ such that for each map
$\delta\colon\Delta[m]\ra\Delta[n]$ of standard simplices, the induced
diagram of pullback squares
$$\xymatrix{
  {P} \ar[r]^h \ar[d]
  & {P'} \ar[r] \ar[d] 
  & {X} \ar[d]^f
  \\
  {S_n\times\Delta[m]} \ar[r]^{1\times\delta}
  & {S_n\times\Delta[n]} \ar[r]
  & {Y}
}$$
is such that $h$ is a weak equivalence of simplicial sheaves.

\item [(5)] For any boolean localization
$p\colon\Sh\bool{B}\ra\topos{E}$, the inverse image $p^\ast f\colon
p^\ast X\ra p^\ast Y$ of $f$ is sharp in $s\Sh\bool{B}$.

\item [(6)] There exists a boolean localization
$p\colon\Sh\bool{B}\ra\topos{E}$ such that the inverse image $p^\ast
f\colon p^\ast X\ra p^\ast Y$ of $f$ is sharp in $s\Sh\bool{B}$.
\end{enumerate}
\end{theorem}

\begin{proof} \forcepar
(1) implies (2):
This follows from (\ref{prop-all-pbs-ho-cart-is-sharp}), and
the fact that $s\topos{E}$ is right proper.

(2) implies (3) and (4):
Let $S_n=Y_n$.

either (3) or (4) implies (5):
This will follow from (\ref{prop-sharp-maps-in-bool-loc}),
since $p^\ast\colon\topos{E}\ra\Sh\bool{B}$ preserves pullbacks and
epimorphisms. 

(5) implies (6):
This is trivial, since every
$\topos{E}$ has a boolean localization
(\ref{subsec-model-cat-simplicial-sheaves}).   

(6) implies (1):
If
$p\colon\Sh\bool{B}\ra\topos{E}$ is a boolean localization, and
$p^\ast f$ is sharp, then since $p^\ast\colon s\topos{E}\ra s\Sh\bool{B}$
preserves pullbacks and reflects weak equivalences
(\ref{subsec-model-cat-simplicial-sheaves}), it follows that $f$ 
is sharp.
\end{proof}

\begin{remark}
In the case when $s\topos{E}=\sSet$, and $f\colon X\ra Y$ a map of
simplicial sets, the above theorem implies that the
following three statements are equivalent.
\begin{enumerate}
\item [(1)]
$f$ is sharp.

\item [(2)]
For each map $g\colon \Delta[n]\ra Y$ the pullback square of $f$ along
$g$ is homotopy cartesian.

\item [(3)]
For each diagram of pullback squares of the form
$$\xymatrix{
  {P} \ar[r]^h \ar[d] 
  & {P'} \ar[r] \ar[d]
  & {X} \ar[d]^f
  \\
  {\Delta[m]} \ar[r]^\delta
  & {\Delta[n]} \ar[r]
  & {Y}
}$$
the map $h$ is a weak equivalence.
\end{enumerate}
Note that characterization (2) is reminiscent of the definition of
quasi-fibration of topological spaces.

A sharp map to a simplicial set $Y$ induces a ``good diagram'' indexed
by the simplices of $Y$, in the sense of
Chach\'olski~\cite{chacholski-closed-classes}.  
\end{remark}

\begin{remark}
Recall from (\ref{exam-boolean-loc}) that if $\topos{E}=\Psh{\cat{C}}$
is a category of presheaves 
on $\cat{C}$, then a suitable boolean localization for $\topos{E}$ is
$\Psh{\cat{C_0}}$.  This implies using part (6) of
(\ref{thm-equiv-charac-of-sharp}) that a map $f\colon X\ra 
Y$ of presheaves on $\cat{C}$ is sharp if and only if for each object
$C\in\cat{C}$ the map $f(C)\colon X(C)\ra Y(C)$ is a sharp map of
simplicial sets.
\end{remark}

\begin{remark}
Recall from (\ref{exam-boolean-loc}) that if $\topos{E}=\Sh(T)$ where
$T$ is a topological space, 
then a boolean localization for $\topos{E}$ is $\Sh T^\delta$.  This
implies using part (6) of (\ref{thm-equiv-charac-of-sharp}) that a map
$f\colon X\ra Y$ of sheaves over $T$ is sharp if and only if for each
point $p\in T$ the induced map $f_p\colon X_p\ra Y_p$ on stalks is a
sharp map of simplical sets.
\end{remark}

\begin{remark}\label{rem-equiv-charac-in-presheaves}
The statement of (\ref{thm-equiv-charac-of-sharp}) remains true if we
replace $s\topos{E}$ with $s\Psh\cat{C}$ equipped with the presheaf
model category structure of (\ref{subsec-presheaf-model-cat}), and
replace boolean localizations $\Sh\bool{B}\ra\topos{E}$ with composite
maps $\Sh\bool{B}\ra\Sh\cat{C}\ra\Psh\cat{C}$.  That this is the case
follows easily from the observation that $f\colon X\ra Y\in
s\Psh\cat{C}$ is sharp if and only if $L^2f\colon L^2X\ra L^2Y\in
s\Sh\cat{C}$ is sharp, the proof of which is straightforward.
\end{remark}

\section{Basic properties of sharp maps of simplicial sheaves}
\label{sec-basic-props-of-sharp-maps}

In this section we give some basic properties of sharp maps in a
simplicial topos.

\begin{theorem}\label{thm-basic-props-of-sharp-maps}
The following hold for simplicial objects in a topos $\topos{E}$.
\begin{enumerate}
\item [P1 ] Local fibrations are sharp.
\item [P2 ] For any object $X\in s\topos{E}$ the projection
map $X\ra 1$ is sharp.
\item [P3 ] Sharp maps are closed under base-change.
\item [P4 ] If $f$ is a map such that
the base-change of $f$ along some epimorphism is sharp, then $f$
is sharp.
\item [P5 ] If maps $f_\alpha$ are sharp for each
$\alpha\in A$ for some set $A$, then the coproduct $\amalg f_\alpha$
is sharp.
\item [P6 ] If $p\colon \topos{E}\ra\topos{E}'$ is a geometric
morphism of topoi, the inverse image functor $p^\ast\colon
s\topos{E}'\ra s\topos{E}$ preserves sharp maps.
\end{enumerate}
\end{theorem}

\begin{proof}
Property P1  follows from part (6) of
(\ref{thm-equiv-charac-of-sharp}), the fact that global fibrations are
sharp (\ref{prop-fibs-are-sharp}), and the fact that local fibrations
are global fibrations in a Boolean localization
(\ref{prop-local-fibs-are-global}).  

Property P2  follows immediately from the fact weak equivalences in
$s\topos{E}$ are precisely those maps $f$ such that $(L^2\Ex^\infty
p^*f)(b)$ is a weak equivalence for each $b\in\bool{B}$ (where $p\colon
\Sh\bool{B}\ra\topos{E}$ is a boolean localization), together with the
fact that the functor $L^2\Ex^\infty p^*$ preserves products.

Property P3  has already been noted in Section~\ref{sec-sharp-maps}.

To prove property P4, consider the pull-back squares
$$\xymatrix{
  {Q} \ar[r]^q \ar[d]
  & {Q'} \ar[r] \ar[d]
  & {P} \ar[r] \ar[d]^g
  & {X} \ar[d]^f
  \\
  {C_n\times\Delta[m]} \ar[r]^{1\times\delta}
  & {C_n\times\Delta[n]} \ar[r]
  & {C} \ar[r]^p
  & {Y}
}$$
where $g$ is sharp and $p$ is an epimorphism.  Then $q$ is a weak
equivalence since $1\times\delta$ is, whence $f$ is sharp by part (4) of
(\ref{thm-equiv-charac-of-sharp}), since the map $C_n\ra Y_n$ is
an epimorphism in $\topos{E}$. 

To prove property P5, let $f_\alpha\colon X_\alpha\ra
Y_\alpha$ be a collection of sharp maps, and let
$f=\amalg_{\alpha\in I}f_\alpha$.  Then P5  follows from the
fact that a coproduct of weak equivalences is a weak equivalence
(\ref{prop-model-cat-colim-lemma}) and
using part (4) of (\ref{thm-equiv-charac-of-sharp}).

Property P6  follows easily from part (4) of
(\ref{thm-equiv-charac-of-sharp}), and the fact that inverse image
functors preserve pullbacks, epimorphisms, and weak equivalences.
\end{proof}

We can now easily prove (\ref{thm-inv-image-ho-cartesian}).
\begin{proof}[Proof of (\ref{thm-inv-image-ho-cartesian})]
Recall (\ref{subsec-model-cat-simplicial-sheaves}) that any homotopy
cartesian square is weakly equivalent to a 
pullback square in which all the maps are fibrations.  Since
fibrations are sharp by (\ref{prop-fibs-are-sharp}), the
square obtained by applying the inverse image functor
$p^*\colon \topos{E}'\ra \topos{E}$ is a pullback square in which the
maps 
are sharp by P6, and hence is a homotopy cartesian square by
(\ref{prop-all-pbs-ho-cart-is-sharp}). 
Since $p^*$ preserves weak equivalences the conclusion follows.
\end{proof}

\begin{remark}
Parts P1--P5 of (\ref{thm-basic-props-of-sharp-maps}) remain true if
we replace $s\topos{E}$ with $s\Psh\cat{C}$ equipped with the presheaf
model category of (\ref{subsec-presheaf-model-cat}), for the reasons
discussed in (\ref{rem-equiv-charac-in-presheaves}).
\end{remark}

\section{Diagonal of a simplicial object}
\label{sec-diagonal-of-s-obj}

Let $X\colon\Delta^{\op}\ra s\topos{E}$ be a simplicial object in
$s\topos{E}$; we write $[n]\mapsto X(n)$ where $X(n)\in s\topos{E}$.
The \dfn{diagonal} 
$\realiz{X}$ of $X$ is an object in 
$s\topos{E}$ defined by $[n]\mapsto X(n)_n$.  

\begin{theorem}\label{thm-sharp-maps-of-simpl-obj}
Let $p\colon X\ra Y$ be a map of simplicial objects in $s\topos{E}$
such that each $p(n)\colon X(n)\ra Y(n)$ is sharp, and each square
$$\xymatrix{
{X(n)} \ar[r] \ar[d]
& {X(m)} \ar[d]
\\ 
{Y(n)} \ar[r]
& {Y(m)}
}$$
is homotopy cartesian.  Then $\realiz{p}\colon
\realiz{X}\ra\realiz{Y}$ is sharp.
\end{theorem}

We prove this theorem using the following well-known inductive
construction of the diagonal of a simplicial object.  Namely,
$\realiz{X}\approx\colim_nF_n\realiz{X}$, where $F_0\realiz{X}=X(0)$ and
each $F_n\realiz{X}$ is obtained from $F_{n-1}\realiz{X}$ by a pushout
diagram of the form
\begin{equation}
\label{eq-pushout-for-diagonal}
\vcenter{\xymatrix{
{X(n)\times\partial\Delta[n]\mskip-20mu
\bigcup_{L_{n-1}X\times\partial\Delta[n]}\mskip-20mu
L_{n-1}X\times\Delta[n]} \ar[r] \ar[d]
& {X(n)\times\Delta[n]} \ar[d]
\\
{F_{n-1}\realiz{X}} \ar[r]
& {F_n\realiz{X}}
}}\end{equation}
where $L_{n-1}X$ denotes the subobject of $X(n)$ which is the union of
the images of all degeneracy maps $s_i\colon X_{n-1}\ra X_n$ for
$0\leq i\leq n$.  

\subsection{Colimits on posets of proper subsets}
\label{subsec-poset-facts}

Before going to the proof of (\ref{thm-sharp-maps-of-simpl-obj}) we
collect some facts about colimits of diagrams indexed by the subsets
of a finite set.  These facts will also be needed in
Sections~\ref{sec-special-diagrams} and
\ref{sec-sharp-maps-in-bool-loc}.  

If $S$ is a finite set, let $\powob{S}$ denote the poset of subsets of $S$,
and let $\propsub{S}$ denote the poset of \emph{proper} subsets
of $S$; we regard $\powob{S}$ and $\propsub{S}$ as categories with $T\ra T'$ if
$T\subseteq T'\subseteq S$.  

Given a functor $X\colon\powob{S}\ra s\topos{E}$ and a subset
$S'\subset S$, we define $X|_{S'}\colon\powob{S'}\ra s\topos{E}$ to be
the restriction of $X$ to $\powob{S'}$ via the formula
$X|_{S'}(T)=X(T)$ for $T\subset S'$.  We also speak of the restriction
$X|_{S'}\colon\propsub{S'}\ra s\topos{E}$ to $\propsub{S'}$.  We say
that a functor 
$X\colon\powob{S}\ra s\topos{E}$ (resp.\ a functor $X\colon \propsub{S}\ra
s\topos{E}$) is \dfn{cofibrant} if for each subset (resp.\ proper
subset) $T\subset S$ the induced map $\colim_{\propsub{T}}X|_T\ra
X(T)$ is a monomorphism.  (The cofibrant functors are in fact the
cofibrant objects in a model category structure on the categories of
functors $\powob{S}\ra s\topos{E}$ and $\propsub{S}\ra s\topos{E}$.)

Say that $S=\{1,\dots,n\}$, and let
$S'=\{1,\dots,n-1\}$.  Define $X'\colon\propsub{S'}\ra
s\topos{E}$ by the formula 
$X'(T)=X(T\cup\{n\})$ for $T\subset S'$.  There is a natural map
$X|_{S'}\ra X'$ of diagrams indexed by $\propsub{S'}$.

\begin{propos}\label{prop-poset-induction}
Suppose $X\colon \propsub{S}\ra s\topos{E}$ is a functor.
There is a natural pushout square
$$\xymatrix{
{\colim_{\propsub{S'}}X|_{S'}} \ar[r] \ar[d]
& {\colim_{\propsub{S'}}X'} \ar[d]
\\
{X(S')} \ar[r]
& {\colim_{\propsub{S}}X}
}$$
and if $f$ is a cofibrant functor, then both vertical
maps in the above square are monomorphisms.
\end{propos}
\begin{proof}
This is a straightforward induction argument on the size of $S$, using
the fact that in a topos, pushouts of monomorphisms are again
monomorphisms. 
\end{proof}

\begin{corol}\label{cor-poset-equivs}
Given a cofibrant functor $X\colon \powob{S}\ra s\topos{E}$ such that
for all $T\subset T'$ the map $X(T)\ra X(T')$ is a weak equivalence, the
induced map $\colim_{\propsub{S}}X\ra X(S)$ is a weak equivalence.
\end{corol}
\begin{proof}
This is proved by a straightforward induction argument on the size of
$S$, using (\ref{prop-poset-induction}).
\end{proof}

\begin{corol}\label{cor-poset-comparison}
Given cofibrant functors $X,Y\colon \propsub{S}\ra s\topos{E}$ and a
map $f\colon 
X\ra Y$ such that each map $X(T)\ra Y(T)$ is a weak equivalence, then
the induced map $\colim_{\propsub{S}}X\ra \colim_{\propsub{T}}Y$ is a
weak equivalence.
\end{corol}
\begin{proof}
This is proved by induction on the size of $S$, using
(\ref{prop-poset-induction}) and the fact that $s\topos{E}$ is a left
proper model category.
\end{proof}

\subsection{The proof of the theorem}

The object $L_{n-1}X$ has an alternate description using the above
notation.  Let
$S=\{1,\dots,n\}$.   Define a functor $F\colon
\powob{S}\ra s\topos{E}$ sending $T\subset S$ to $X(\#T)$,
and sending $i\colon T\ra T'$ to the map induced by the simplicial operator
$\sigma\colon[\#T']\ra[\#T]$ defined by $\sigma(0)=0$ and
$\sigma(k)=\max(\ell\ |\ i(k)\leq\ell)$ for $0<k\leq \#T'$.  Then
$L_{n-1}X=\colim_{\propsub{S}}F$.  Furthermore, $F$ is a cofibrant functor

\begin{proof}[Proof of (\ref{thm-sharp-maps-of-simpl-obj})]
Each map $\partial\Delta[n]\ra\Delta[n]$ is mono, as are the maps
$L_{n-1}X\ra X(n)$, and the top horizontal arrow in
(\ref{eq-pushout-for-diagonal}).  The proof is a
straightforward induction following the inductive 
construction of diagonal given above and using
(\ref{prop-special-diagrams}) together with
(\ref{prop-poset-induction}).  

That is, suppose by induction that
$F_{n-1}\realiz{X}\ra F_{n-1}\realiz{Y}$ is sharp.  Using
(\ref{prop-special-diagrams}, 3) one shows that $L_{n-1}X\ra L_{n-1}Y$
is sharp.  Then using (\ref{prop-special-diagrams}, 2) one shows that
the induced map from the upper left-hand corner of
(\ref{eq-pushout-for-diagonal}) to the upper left-hand corner of the
corresponding square for $Y$ is sharp.  Applying
(\ref{prop-special-diagrams}, 2) to the whole square
(\ref{eq-pushout-for-diagonal}) gives that $F_n\realiz{X}\ra
F_n\realiz{Y}$ is sharp.  
Finally, (\ref{prop-special-diagrams}, 1)
shows that $\realiz{X}\ra \realiz{Y}$ is sharp, as desired.
\end{proof}

\begin{remark}\label{rem-ho-invariance-of-diag}
If $f\colon X\ra Y$ is a map of simplicial objects in $s\topos{E}$ such
that in each degree $n$ the map $f(n)\colon X(n)\ra Y(n)$ is a weak
equivalence, then one may show by using the above inductive scheme
together with (\ref{cor-poset-comparison}) that
$\realiz{f}\colon\realiz{X}\ra\realiz{Y}$ is a weak equivalence, since
$s\topos{E}$ is a proper model category and the cofibrations are
precisely the monomorphisms.
\end{remark}

\section{Homotopy colimits}
\label{sec-ho-colims}

Let $X\colon\cat{I}\ra s\topos{E}$ be a diagram of simplicial
sheaves.  As in \cite{yellow-monster} the \dfn{homotopy colimit} of
$X$, denoted 
$\hocolim_{\cat{I}}X$, is defined to be the 
diagonal of the simplicial object in $s\topos{E}$ given in each degree
$n\geq0$ by
$$[n]\mapsto \coprod_{i_0\ra\cdots\ra i_n}\mskip-20mu Xi_0,$$
where the coproduct is taken over all composable strings of arrows in
$\cat{I}$ of length $n$.
From (\ref{prop-model-cat-colim-lemma}) and
(\ref{rem-ho-invariance-of-diag}) it follows that 
$\hocolim_{\cat{I}}X$ is a weak homotopy equivalence invariant of $X$.

Let $\overcat{\cat{I}}{i}$ denote the category of objects over a
fixed object $i$ in $\cat{I}$.  Given an $\cat{I}$-diagram $X$, one
can define an $\cat{I}$-diagram $\widetilde{X}$ by
$\widetilde{X}i=\hocolim_{\overcat{\cat{I}}{i}}X$.
Thus, $\widetilde{X}i$ is the diagonal of the simplicial object in
$s\topos{E}$ given by
$$[n]\mapsto \coprod_{i_0\ra\cdots\ra i_n\ra i}\mskip-20mu Xi_0.$$
There is a natural map $\widetilde{X}\ra X$ of $\cat{I}$-diagrams, and
an isomorphism of simplicial sheaves
$\hocolim_{\cat{I}}X\approx\colim_{\cat{I}}\widetilde{X}$.  (This is
the construction of \cite{yellow-monster}.)

\begin{theorem}
\label{thm-hocolim-of-equifibered-sharp-is-sharp}
Let $f\colon X\ra Y$ be a map of $\cat{I}$-diagrams of simplicial
sheaves such that
\begin{enumerate}
\item [(1)]
each map $fi\colon Xi\ra Yi$ is sharp for $i\in\ob\cat{I}$, and
\item [(2)]
each square $\vcenter{\xymatrix{
{Xi} \ar[r] \ar[d]
& {Xj} \ar[d]
\\
{Yi} \ar[r]
& {Yj}
}}$ for $\alpha\colon i\ra j\in\cat{I}$ is homotopy cartesian.
\end{enumerate}
Then
\begin{enumerate}
\item [(a)]
the induced map $\hocolim_{\cat{I}}f\colon\hocolim_{\cat{I}}X
\ra\hocolim_{\cat{I}}Y$ is sharp, and
\item [(b)]
for each object $i$ in $\cat{I}$ the square $\vcenter{\xymatrix{
{\widetilde{X}i} \ar[r] \ar[d]
& {\hocolim X} \ar[d]
\\ 
{\widetilde{Y}i} \ar[r] 
& {\hocolim Y}
}}$ is a pull-back square.
\end{enumerate}
\end{theorem}
\begin{proof}
First, we note that (b) follows without need of the hypotheses (1) and
(2).  This is because for each $n\geq 0$, the square
$$\xymatrix{
{\coprod_{i_0\ra\cdots\ra i_n\ra i}\mskip-20mu Xi_0}\ar[r]\ar[d]
&{\coprod_{i_0\ra\cdots\ra i_n}\mskip-20mu Xi_0}\ar[d]
\\
{\coprod_{i_0\ra\cdots\ra i_n\ra i}\mskip-20mu Yi_0}\ar[r]
&{\coprod_{i_0\ra\cdots\ra i_n}\mskip-20mu Yi_0}
}$$
is a pullback square by the distributive law
(\ref{prop-distributive-law}), and taking diagonals of bi-simplicial
objects 
commutes with limits.

To show (a), we consider the square
$$\xymatrix{
{\coprod_{i_0\ra\cdots\ra i_n}\mskip-20mu Xi_0}\ar[r]\ar[d]
&{\coprod_{j_0\ra\cdots\ra j_m}\mskip-20mu Xj_0}\ar[d]
\\
{\coprod_{i_0\ra\cdots\ra i_n}\mskip-20mu Yi_0}\ar[r]
&{\coprod_{j_0\ra\cdots\ra j_m}\mskip-20mu Yj_0}
}$$
where the horizontal arrows are induced by a map $\delta\colon
[m]\ra[n]\in\Delta$.  The vertical arrows are sharp by
(\ref{thm-basic-props-of-sharp-maps}, P5),
and the square is homotopy cartesian using
(\ref{lemma-coproduct-of-ho-cartesian-squares}) .
(In fact, the square is a \emph{pullback} square except when $\delta$
is a simplicial operator for which $\delta(0)\neq0$, in which case the
square is only 
homotopy cartesian.)  The result then follows from
(\ref{thm-sharp-maps-of-simpl-obj}). 
\end{proof}

\begin{lemma}\label{lemma-coproduct-of-ho-cartesian-squares}
In $s\topos{E}$, an arbitrary coproduct of homotopy cartesian squares
is homotopy cartesian.
\end{lemma}
\begin{proof}
A coproduct of weak equivalences is a weak equivalence by
(\ref{prop-model-cat-colim-lemma}),  
and a coproduct of pullback squares is a pullback square by the
distributive law (\ref{prop-distributive-law}).  Thus it suffices to
factor the sides of each square into a weak equivalence followed by a
fibration and demonstrate the result for the resulting pullback
squares; since fibrations are sharp (\ref{prop-fibs-are-sharp}), 
the coproduct of sharp maps is sharp
(\ref{thm-basic-props-of-sharp-maps}, P5), and pullbacks along sharp
maps are homotopy cartesian, the result follows.
\end{proof}

\begin{proof}[Proof of (\ref{thm-hocolims})]
To prove (1), 
choose a factorization
$$\colim_{\cat{I}}X\xra{j}W'\xra{p}\colim_{\cat{I}}Y$$ 
such that $j$ is
a weak equivalence and $p$ is sharp (e.g., a fibration).  Define an
$\cat{I}$-diagram $X'$ 
by $X'i=W'\!\!\pullback{\colim Y}\!\!Yi$; by the distributive law
(\ref{prop-distributive-law}) we
see that 
$\colim_{\cat{I}}X'\approx W'$.  Note also that the induced map $Xi\ra
X'i$ is a weak equivalence, since $p$ is sharp and by the hypothesis
that each square (\ref{eq-square-of-maps-to-colims}) is homotopy
cartesian. 
In the diagram
$$\xymatrix{
{\colim \widetilde{X}} \ar[r]_\sim^i \ar[d]^k
& {\colim \widetilde{X'}} \ar[r] \ar[d]
& {\colim \widetilde{Y}} \ar[d]^\sim_\ell
\\
{\colim X}\ar[r]_{\sim}^j
&{\colim X'}\ar[r]^p
& {\colim Y}
}$$
the map $p$ is sharp and the indicated maps are weak equivalences;
that $i$ and $\ell$ are weak equivalences follows from the homotopy
invariance of homotopy colimits and  the hypothesis that $Y$
is a homotopy colimit diagram.
Thus to show that $k$ is a weak equivalence, and hence that $X$ is a homotopy
colimit diagram, it suffices 
to show that 
the right-hand square is a pull-back square.  

Since each $X'i$ is defined to be the pullback of $\colim X'\ra\colim
Y$ along a map 
$Yi\ra\colim Y$, we see that $\widetilde{X}i$ is the pullback of
$\colim X'$
along the composite map $\widetilde{Y}i\ra Yi\ra \colim Y$.  The
assertion that the right-hand square is a pullback square now follows
using the distributive law (\ref{prop-distributive-law}).

To prove (2), choose a factorization of $f\colon X\ra Y$ into
$X\xra{j}X'\xra{p}Y$, in which $j$ is an object-wise weak equivalence 
and $p$ is an object-wise fibration, and hence object-wise sharp.
Then the square 
$$\xymatrix{ 
{\widetilde{X'}i} \ar[r] \ar[d]
&{\colim_{\cat{I}}\widetilde{X'}}\ar[d]
\\
{\widetilde{Y}i} \ar[r]
&{\colim_{\cat{I}}\widetilde{Y}}
}$$
is homotopy cartesian by
(\ref{thm-hocolim-of-equifibered-sharp-is-sharp}), and since $X$ is by
hypothesis a homotopy colimit diagram it follows that this
square is weakly equivalent to
(\ref{eq-square-of-maps-to-colims}), and we get the desired result.
\end{proof}

\section{Lemmas on sharp maps of special diagrams}
\label{sec-special-diagrams}

In this section we show that for special kinds of maps of diagrams,
the induced map of colimits is sharp.  These results were the key
element of 
the proof of (\ref{thm-sharp-maps-of-simpl-obj}).

\begin{propos}\label{prop-special-diagrams}
Let $\cat{I}$ denote a small category, and $p\colon X\ra Y$ a map of
$\cat{I}$-diagrams in $s\topos{E}$.  Suppose that $pi\colon Xi\ra
Yi$ is sharp for each $i\in\ob\cat{I}$, and that
$$\xymatrix{
{Xi} \ar[r]^{X\alpha} \ar[d]_{pi}
& {Xj} \ar[d]^{pj}
\\
{Yi} \ar[r]^{Y\alpha}
& {Yj}
}$$
 is a homotopy cartesian square for each $\alpha\colon i\ra j$ in
$\cat{I}$.  Then under each of the following cases (1)--(3), the
induced map $\colim_{\cat{I}}X\ra\colim_{\cat{I}}Y$ is sharp, and for
each $i\in\ob\cat{I}$, the square
$$\xymatrix{
{Xi} \ar[r] \ar[d]
& {\colim_{\cat{I}}X} \ar[d]
\\
{Yi} \ar[r]
& {\colim_{\cat{I}}Y}
}$$ 
is homotopy cartesian.
\begin{enumerate}
\item [(1)] $\cat{I}$ is the category obtained from the poset $\N$ of
non-negative integers, and each map $X(n)\ra X(n+1)$ and $Y(n)\ra
Y(n+1)$ is a monomorphism.

\item [(2)] $\cat{I}=(i_1\xleftarrow{\alpha}i_0\xra{\beta}i_2)$ and
$X\beta$ and $Y\beta$ are monomorphisms.

\item [(3)] $\cat{I}$ is the category obtained from the poset
$\bar{\mathcal{P}}S$ of
\emph{proper} subsets of a finite set $S$, and $X$ and $Y$ are
cofibrant diagrams in the sense of (\ref{subsec-poset-facts}).
\end{enumerate}
\end{propos}

%

\begin{lemma}\label{lemma-silly-sequence-arg}
Consider a countable sequence of maps over $B$
$$Y(0)\xra{i_0}Y(1)\xra{i_1}Y(2)\xra{i_2}\dots\ra B$$
such that each $i_n$ is a trivial cofibration, and each map $q_n\colon
Y(n)\ra B$ is sharp.  Then the induced map $q\colon\colim_n Y(n)\ra B$
is sharp. 
\end{lemma}
\begin{proof}
Given a map $f\colon A\ra B$, consider the pullbacks $X(n)=Y(n)\times_{B}A$.
By the distributive law (\ref{prop-distributive-law}),
$\colim_nX(n)=\colim_nY(n)\times_{B}A$.  Since each map $Y(n)\ra B$ is
sharp and each $i_n$ is a weak equivalence, it follows that each
$X(n)\ra X({n+1})$ is a weak equivalence 
and thus a trivial cofibration.
Thus the composite $X(0)\ra \colim_nX(n)$ is a trivial cofibration,
and so base-change  
of $q$ along $f$ yields a homotopy cartesian square.
\end{proof}

\begin{proof}[Proof of part 1 of (\ref{prop-special-diagrams})]
Let $X'(n)=Y(n)\!\!\times_{\colim Y(n)}\!\!\colim X(n)$, whence we
have that $\colim
X'(n)\approx\colim X(n)$ by the distributive law
(\ref{prop-distributive-law}).  It suffices to show 
\begin{enumerate}
\item [(1)]
that each map
$X'(n)\ra Y(n)$ is sharp, and
\item [(2)]
that each map
$X(n)\ra X'(n)$ is a weak equivalence.
\end{enumerate}
This is because (1), together with (\ref{thm-equiv-charac-of-sharp},
P4) and the fact that $\coprod_n Y(n)\ra\colim Y(n)$ is epi, implies
that $\colim X(n)\ra \colim Y(n)$ is sharp, and (2) then demonstrates
that the appropriate squares are homotopy cartesian. 

Let $X(n,m)=Y(n)\times_{Y(m)}X(m)$ for $m\geq n$.  Then
$X'(n)\approx\colim_m 
X({n,m})$ by the distributive law.  We have that $X(n,n)\approx
X(n)$, and 
each map $X(n,m)\ra X(n,m+1)$ is a weak equivalence since $X(m+1)\ra
Y(m+1)$ is sharp.  Thus $X(n)\ra \colim_m X(n,m)\approx X'(n)$ is a weak
equivalence, proving (2).
Claim (1) follows from
(\ref{lemma-silly-sequence-arg}) applied to the sequence $X({n,m})$
over $Y(n)$.
\end{proof}

The following lemma describes conditions under which one may ``glue''
an object onto a sharp map and still obtain a sharp map.
\begin{lemma}\label{lemma-silly-pushout-arg}
Let
$$\xymatrix{
  {X} \ar[r]^i \ar[d]_f
  & {Y} \ar[d]^g
  \\
  {X'} \ar[r] \ar[d]_p
  & {Y'} \ar[d]^q
  \\
  {A} \ar[r]
  & {B}
}
$$
be a commutative diagram such that the top square is a push-out
square, $p$, $pf$, and $qg$ 
are sharp, $f$ is a weak equivalence, and either $i$ or $f$ is a
monomorphism.  Then $q$ is also sharp. 
\end{lemma}
\begin{proof}
It suffices by (\ref{prop-all-pbs-ho-cart-is-sharp}) to show
that every base-change of $q$ along a map $U\ra B$ 
produces a homotopy cartesian square.  Since $qg$ is sharp it
suffices to show that 
$$U\pullback{B}g\colon U\pullback{B}Y\ra U\pullback{B}Y'$$ 
is a weak equivalence.  Via the pushout square
$$\xymatrix{
{U\pullback{B}{X}} \ar[r]_{U\pullback{B}i} \ar[d]_{U\pullback{B}f}
& {U\pullback{B}{Y}} \ar[d]
\\
{U\pullback{B}{X'}} \ar[r]
& {U\pullback{B}{Y'}}
}$$
in which either the top or the left arrow is a cofibration, we see
that it
suffices to show that $U\times_{B}f$ is a weak equivalence, since
$s\topos{E}$ is a left-proper model category.

In fact, $U\times_{B}f\approx (U\times_{B}A)\times_{A}f$; that
is, $U\times_Bf$  is a base-change of $f$ along a map into $A$.  Thus
since $p$ and $pf$ are sharp, this base-change of $f$ is a weak
equivalence, as desired.
\end{proof}

We have need of the following peculiar lemma.
\begin{lemma}\label{lemma-peculiar}
In a Grothendieck topos $\topos{E}$ consider a diagram of the form
$$\xymatrix{
{A} \ar@{>->}[r] \ar[d]^p
& {X} \ar[d]
\\ 
{A'} \ar@{>->}[r] \ar[d]
& {X'} \ar[d]^q
\\
{B} \ar@{>->}[r] 
& {Y}
}$$
in which the horizontal arrows are mono, the top square is a pushout
square, and the large rectangle is a pullback rectangle.  Then the bottom 
square is also a pullback square.
\end{lemma}
\begin{proof}
It suffices to show that the lemma holds in a Boolean localization of
$\topos{E}$.  In this case every subobject has a complement, so we may
write $X=A\amalg C$, $X'=A'\amalg C'$, and $Y=B\amalg D$.  To show
that the lower square is a pullback, it suffices to show that
$q(C')\subset D$.  Since the top square is a pushout, $p(C)=C'$, and
since the big rectangle is a pullback, $qp(C)\subset D$, producing the
desired result.
\end{proof}

\begin{proof}[Proof of part 2 of (\ref{prop-special-diagrams})]
We have a diagram of the form
\begin{equation}\label{eq-map-of-pushouts}
\vcenter{\xymatrix{
{X_1} \ar[d]_{p_1}
& {X_0} \ar[l] \ar[d]_{p_0} \ar@{>->}[r]^{i}
& {X_2} \ar[d]^{p_2}
\\
{Y_1} 
& {Y_0} \ar[l] \ar@{>->}[r]^{j}
& {Y_2}
}}
\end{equation}
where $p_n$ is sharp for $n=0,1,2$, each square is homotopy cartesian,
and $i$ and $j$ are mono.  We must show that the induced map
$X_{12}\ra Y_{12}$ of pushouts is sharp, and that each square
\begin{equation}\label{eq-pullback-from-pushout}
\vcenter{\xymatrix{
{X_n} \ar[r] \ar[d]_{p_n}
& {X_{12}} \ar[d]^{p_{12}}
\\
{Y_n} \ar[r]
& {Y_{12}}
}}
\end{equation} is homotopy cartesian for $n=0,1,2$.

We prove the claim by proving it for the following cases:
\begin{enumerate}
\item [(a)] under the additional hypothesis that both of the squares in
(\ref{eq-map-of-pushouts}) are pullback squares, 
\item [(b)] under the additional hypothesis that the right-hand square
in (\ref{eq-map-of-pushouts}) is a pull back square, and
\item [(c)] under no additional hypotheses.
\end{enumerate}

In case (a), each square of the form (\ref{eq-pullback-from-pushout})
is necessarily a pullback square since $i$ and $j$ are mono; this can
be seen by passing to a boolean localization, in which case
$X_2\approx X_0\amalg X_2'$ and $Y_2\approx Y_0\amalg Y_2'$ so that
$X_{12}\approx X_1\amalg X_2'$ and $Y_{12}\approx Y_1\amalg Y_2'$.
Thus $p_{12}\colon X_{12}\ra Y_{12}$ must be sharp using
(\ref{thm-equiv-charac-of-sharp}, P4), since  the pullback
of $p_{12}$ along the epimorphism $Y_1\amalg Y_2\ra Y_{12}$ is sharp.

In case (b), we let $X_0'=Y_0\pullback{Y_1}X_1$ and
$X_2'=X_0'\union{X_0}X_2$, obtaining a diagram of the form
$$\xymatrix{
& {X_0} \ar@{>->}[r]^i \ar[d]_{\sim}
& {X_2} \ar[d]_\sim
\\
{X_1} \ar[d]_{p_1}
& {X_0'} \ar[l] \ar@{>->}[r] \ar[d]_{p_0'}
& {X_2'} \ar[d]^{p_2'}
\\
{Y_1} 
& {Y_0} \ar[l] \ar@{>->}[r]^j
& {Y_2}
}$$
The map $p_0'$ is sharp since it is a
base-change of the sharp map $p_1$.  
The 
map $p_2'$ is sharp by (\ref{lemma-silly-pushout-arg}) since $i$
is mono.  
The lower right-hand square is a pullback square by
(\ref{lemma-peculiar}).
Then the claim
reduces to case 
(a), since $X_1\cup_{X_0'}X_2'\approx X_{12}$.  

In case (c), let $X_0'=Y_0\pullback{Y_2}X_2$ and
$X_1'=X_1\union{X_0}X_0'$, obtaining a diagram of the form
$$\xymatrix{
{X_1} \ar@{>->}[d]_\sim
& {X_0} \ar[l] \ar@{>->}[d]^{\sim}
\\
{X_1'} \ar[d]_{p_1'} 
& {X_0'} \ar[l] \ar[d]_{p_0'} \ar@{>->}[r]^{i'}
& {X_2} \ar[d]^{p_2}
\\
{Y_1} 
& {Y_0} \ar[l] \ar@{>->}[r]^j
& {Y_2}
}$$
The map $p_0'$ is the base-change of a sharp
map $p_2$ and hence is sharp; the map $p_1'$ is sharp by (\ref{lemma-silly-pushout-arg}) (note
that $X_0\ra X_0'$ is mono).  Thus the claim reduces to case (b), since
$X_1'\cup_{X_0'}X_2\approx X_{12}$.
\end{proof}

\begin{proof}[Proof of part 3 of (\ref{prop-special-diagrams})]
Let $S=\{1,\dots,n\}$; we prove the result by induction on $n$.  The
cases $n=0,1$ are trivial, and case $n=2$ follows from
(\ref{prop-special-diagrams}, part 2).  

For a set $T$, as in (\ref{subsec-poset-facts}) let $\bar{P}T$ denote
the poset of \emph{proper} 
subsets of $T$ as in (\ref{subsec-poset-facts}).  Then
(\ref{prop-poset-induction}) provides a pushout square 
$$\xymatrix{
{\colim_{\bar{P}S'}X|_{S'}} \ar[r] \ar@{>->}[d]
& {\colim_{\bar{P}S'}X'} \ar@{>->}[d]
\\
{X(S')} \ar[r]
& {\colim_{\bar{P}S} X}
}$$
in which the vertical arrows are mono; here $S'=\{1,\dots,n-1\}$.
There is a similar 
diagram for $Y$.  One now deduces the result by induction on the size
of $S$, applying (\ref{prop-special-diagrams}, part 2) to the above
square to carry out the induction step.

Note that in order to apply the induction step, we need to know that
the square 
$$\xymatrix{
{\colim_{\propsub{S'}}X|_{S'}} \ar[r] \ar[d]
& {\colim_{\propsub{S'}}X'} \ar[d]
\\
{\colim_{\propsub{S'}}Y|_{S'}} \ar[r]
& {\colim_{\propsub{S'}}Y'}
}$$
is homotopy cartesian.  This follows by induction from the fact that
each of the squares in
$$\xymatrix{
{\colim_{\propsub{S'}}X|_{S'}} \ar[r] \ar[d]
& {X(S')} \ar[r] \ar[d]
& {X(S)} \ar[d]
& {\colim_{\propsub{S'}}X'} \ar[l] \ar[d]
\\
{\colim_{\propsub{S'}}Y|_{S'}} \ar[r]
& {Y(S')} \ar[r]
& {Y(S)}
& {\colim_{\propsub{S'}}Y'} \ar[l]
}$$
are homotopy cartesian.
\end{proof}

\section{Sharp maps in a boolean localization}
\label{sec-sharp-maps-in-bool-loc}

In this section we go back to prove the results needed in the proof of
(\ref{thm-equiv-charac-of-sharp}) on sharp maps in a boolean
localization. 
In the following $\bool{B}$ denotes a complete boolean algebra.

\begin{propos}
\label{prop-sharp-maps-in-bool-loc}
Let $f\colon X\ra Y$ be a map in $s\Sh\bool{B}$.  The following are
equivalent.
\begin{enumerate}
\item [(1)]
$f$ is sharp.

\item [(2)]
For all $n\geq 0$ and all $S_n\ra Y_n$ in $\Sh\bool{B}$ the induced
pullback square
$$\xymatrix{
  {P} \ar[r] \ar[d]
  & {X} \ar[d]^f
  \\
  {S_n\times\Delta[n]} \ar[r]
  & {Y}
}$$
is homotopy cartesian.

\item [(3)]
For each $n\geq 0$ there exists an epimorphism $S_n\ra Y_n$ in
$\Sh\bool{B}$ such that for each map $\delta\colon
\Delta[m]\ra\Delta[n]$ of standard simplices, the induced diagram of
pullback squares  
$$\xymatrix{
  {P} \ar[r]^h \ar[d]
  & {P'} \ar[r] \ar[d]
  & {X} \ar[d]^f
  \\
  {S_n\times\Delta[m]} \ar[r]^{1\times\delta}
  & {S_n\times\Delta[n]} \ar[r]
  & {Y}
}$$
is such that $h$ is a weak equivalence.
\end{enumerate}
\end{propos}

Let $\gamma$ be an ordinal, viewed as a category.  Given a
functor $X\colon \gamma\ra s\topos{E}$ such that
$\colim_{\alpha<\beta}X(\alpha)\approx X(\beta)$ for each limit
ordinal $\beta<\gamma$, we call the induced map $X(0)\ra
\colim_{\alpha<\gamma}X(\alpha)$ a \dfn{transfinite composition} of
the maps in 
$X$. 

\begin{proof}
The implications (1) implies (2) implies (3)  are
straightforward, so it suffices to prove (3) implies (1).

Consider a diagram of pullback squares
\begin{equation}
\label{eq-another-double-pullback}
\vcenter{\xymatrix{
  {A} \ar[r]^g \ar[d]
  & {A'} \ar[r] \ar[d]
  & {X} \ar[d]^f
  \\
  {B} \ar[r]^h 
  & {B'} \ar[r]
  & {Y}
}}
\end{equation}
with $h$ a weak equivalence and $f$ as in (3).  We want to show that $g$
is a weak 
equivalence.  By factoring $h$ into a trivial cofibration followed by
a trivial fibration, we see that we can reduce to the case when $h$ is
a trivial cofibration. 

Let $\mathcal{C}$ denote the class of trivial cofibrations $h\colon
B\ra B'$ such that for all diagrams of the form
(\ref{eq-another-double-pullback}) in which the squares are pullbacks,
the map $g$ is a weak equivalence.  In order to show that
$\mathcal{C}$ contains \emph{all} trivial cofibrations, it will
suffice by (\ref{lemma-classc-is-tcof}) to show that
$\mathcal{C}$ 
\begin{enumerate}
\item [(1)] is closed under retracts,
\item [(2)] is closed under cobase-change,
\item [(3)] is closed under transfinite composition, and
\item [(4)] contains all maps of the form $U\times\Lambda^k[n]\ra
U\times\Delta[n]$ where $U\in\topos{E}$ is a discrete object.
\end{enumerate}

Part (1) is straightforward.  Part (2) follows from the fact that
pullbacks of monomorphisms are monomorphisms, and the fact that the
cobase-change of a trivial cofibration is again a trivial
cofibration.  Part (3) follows from the fact that a transfinite
composite of trivial cofibrations is a trivial cofibration, and from
the fact that
if $\{Y_\alpha\}$ is some sequence and $f\colon X\ra Y=\colim_\alpha
Y_\alpha$ a map, then 
$\colim_\alpha Y_\alpha\times_{Y}X\approx X$ by the distributive law
(\ref{prop-distributive-law}). 

Part (4) is (\ref{lemma-horn-gluing}).
\end{proof}





Let $\Lambda^k[n]\subset\Delta[n]$ denote the ``$k$-th horn'' of the
standard $n$-simplex; that is, $\Lambda^k[n]$ is the largest
subcomplex of $\Delta[n]$ not containing the $k$-th face.  The
following lemma, though simple, is crucial to proving anything about
sharp maps.  It is essentially Lemma~7.4 of Chach\'olski
\cite{chacholski-closed-classes}, at least in the case when
$\Sh\bool{B}=\Set$. 

\begin{lemma}
\label{lemma-horn-gluing}
Let $f\colon X\ra Y$ be a map in $s\Sh\bool{B}$, and let $U_n\ra
Y_n$ be some map in $\Sh\bool{B}$.  Suppose that for each map
$\delta\colon \Delta[m]\ra\Delta[n]$ of standard simplices the map $i$
in the diagram
$$\xymatrix{
  {P_\delta} \ar[r]^i \ar[d]
  & {P} \ar[r] \ar[d]
  & {X} \ar[d]^f
  \\ 
  {U_n\times\Delta[m]} \ar[r]^{1\times\delta}
  & {U_n\times\Delta[n]} \ar[r]
  & {Y}
}$$
of pullback squares
is a weak equivalence.  Then for any inclusion
$\Lambda^k[n]\ra\Delta[n]$ of a simplicial horn into a standard
simplex the map $j$ in the diagram
$$\xymatrix{
  {Q} \ar[r]^j \ar[d]
  & {P} \ar[r] \ar[d]
  & {X} \ar[d]^f
  \\ 
  {U_n\times\Lambda^k[n]} \ar[r]^{1\times\delta}
  & {U_n\times\Delta[n]} \ar[r]
  & {Y}
}$$
of pullback squares is a weak equivalence.
\end{lemma}
\begin{proof}
Let $S=\{0,\dots,n\}$ be a set, and identify $S$ with the set of
vertices of $\Delta[n]$.  There is a functor $F\colon
\powob{S}\ra\sSet$ sending $T\subseteq S$ to the smallest 
subobject of $\Delta[n]$ containing the vertices $T$; thus
$F(T)\approx\Delta[\realiz{T}]$.  Note that $F$
is a cofibrant functor in the sense of (\ref{subsec-poset-facts}).
Let
$S'=S\setminus\{k\}$, and define
$\widetilde{F}\colon \powob{S'}\ra \sSet$ by
$\widetilde{F}(T)=F(T\cup\{k\})$.  Then
$\colim_{\propsub{S'}}\widetilde{F}\approx\Lambda^k[n]$, and
$\widetilde{F}(S')\approx\Delta[n]$, and $\widetilde{F}$ is a 
cofibrant functor.

Now define $G\colon\powob{S'}\ra s\Sh\bool{B}$ by $G(T)=X\times_Y
(U_n\times \widetilde{F}(T))$.  Then $\colim_{\propsub{S'}}G\approx Q$
by the distributive law (\ref{prop-distributive-law})
and $G(S')\approx P$, and the lemma follows immediately from
(\ref{cor-poset-equivs}), since $G$ is a cofibrant functor.
\end{proof}

%
%
%
%

\section{Local fibrations are global fibrations in a boolean
  localization}
\label{sec-local-fibs-are-global-fibs}

The purpose of this section is to prove
(\ref{prop-local-fibs-are-global}), as well as
(\ref{lemma-classc-is-tcof}), which was used in
Section~\ref{sec-sharp-maps-in-bool-loc}.   
It is possible with some work to derive these facts from Jardine's
construction of the model category on $s\Sh\bool{B}$.  However, it
seems more enlightening (and no more difficult) to proceed by
constructing the model category 
structure on $s\Sh\bool{B}$ from scratch, and showing that it has the
desired properties while coinciding with Jardine's structure; it turns
out that the construction of the model category structure on
$s\Sh\bool{B}$ is somewhat simpler than the more general case of
simplicial sheaves on an arbitrary Grothendieck site.

\subsection{Sheaves on a complete boolean algebra}

Let $\bool{B}$ be a \dfn{complete boolean algebra}.  Thus $\bool{B}$ is a
complete distributive lattice with minimal and maximal elements $0$
and $1$, such 
that every $b\in \bool{B}$ has a complement $\bar{b}$; that is, if $\vee$
denotes meet and $\eev$ denotes join, then $b\vee \bar{b}=1$ and $b\eev
\bar{b}=0$.  We view $\bool{B}$ as a category, with a map $b\ra b'$
whenever $b\leq b'$.  

A \dfn{presheaf} on $\bool{B}$ is a functor $X\colon \bool{B}^\op\ra
\Set$.  A \dfn{sheaf} on $\bool{B}$ is a presheaf $X$ such that for
a collection of elements $\{b_i\in \bool{B}\}_{i\in I}$, the diagram
\begin{equation}\label{eq-sheaf-equalizer}
X(b)\lra \prod_{i\in I}X(b_i) \rightrightarrows\prod_{i,j\in
I}X(b_i\eev b_j)
\end{equation}
is an equalizer whenever $\bigvee_{i\in I}b_i=b$.

Say a collection of elements $\{b_i\in \bool{B}\}_{i\in I}$ is a
\dfn{decomposition} of $b\in\bool{B}$ if $\bigvee_{i\in I}b_i=b$ and $b_i\eev
b_j=0$ if $i\neq j$.  We write $\coprod_{i\in
I}b_i=b$ to denote a 
decomposition of $b$.  The collection of decompositions of
$b\in\bool{B}$ forms a 
directed set under refinement, with the trivial decomposition $\{b\}$
as minimal element.

\begin{propos}\label{prop-bool-loc-sheaf-char}
A presheaf $X\in\Psh\bool{B}$ is a sheaf if and only if for each
decomposition $b=\coprod_{i\in I} b_i$ of $b$, the induced map
$X(b)\ra\prod_{i\in I}X(b_i)$ is an isomorphism.
\end{propos}
\begin{proof}
The only if statement follows from the definition of a sheaf.  To
prove the if statement, suppose 
$\bigvee_{i\in I}b_i=b$.  For each $S\subset I$, let
$$b_S=b\eev\left(\bigwedge_{i\in S}b_i\right)
\eev\left(\bigwedge_{j\in S\setminus I}\bar{b}_j\right).$$
Note that for $T\subset I$ we have that $\bigvee_{i\in
T}b_i=\coprod_{S\subset T}b_S$ and $\bigwedge_{i\in
T}b_i=\coprod_{S\supset T}b_S$.  In particular, $b_\varnothing=0$ and
$b=\coprod_{S\subset I}b_S$.

To show that $X$ is a sheaf, we need to show that every sequence of
the form (\ref{eq-sheaf-equalizer})
is exact.  By hypothesis this sequence is isomorphic to the diagram
$$\prod_{S\subset I}X(b_S)\ra \prod_{i\in I}\prod_{T\ni i}X(b_T)
\rightrightarrows\prod_{i,j\in I}\prod_{U\ni i,j}X(b_U).$$
But the above sequence is manifestly exact, since it is a product of
sequences of the form
$$X(b_S)\ra\prod_{i\in S}X(b_S) \rightrightarrows\prod_{i,j\in
S}X(b_S)$$
for each $S\subset I$.
\end{proof}

Define a functor $L\colon \Psh\bool{B}\ra \Psh\bool{B}$ by
$$(LX)(b)=\coliml_{b=\amalg b_i} \prod_i X(b_i),$$
the colimit being taken over the directed set of decompositions of
$b$.  There is a natural transformation $\eta\colon X\ra LX$
corresponding to the trivial decompositions of $b\in\bool{B}$.
Typically, one proves that the composite functor $L^2=L\circ L$
is a sheafification functor for sheaves on $\bool{B}$.  The following
shows that $L$ \emph{is itself a sheafification functor.}

\begin{propos}\label{prop-sheafification-bool-loc}
Given $X\in\Psh\bool{B}$, the object $LX$ is a sheaf; furthermore,
$\eta_X\colon X\ra LX$ is an isomorphism if $X$ is a sheaf.  Thus $L$
induces a sheafification functor
$L\colon\Psh\bool{B}\ra\Sh\bool{B}$ left adjoint to inclusion
$\Sh\bool{B}\ra\Psh\bool{B}$. 
\end{propos}
\begin{proof}
The proof is a straightforward element chase, using
(\ref{prop-bool-loc-sheaf-char}) and the fact that decompositions of
an element $b\in\bool{B}$ form a directed set.
%
%
\end{proof}

\subsection{A model category structure for $s\Sh\bool{B}$}

Let $f\colon X\ra Y$ be a map in $s\Sh\bool{B}$.  Say that $f$ is 
\begin{enumerate}
\item a \dfn{cofibration} if it is a monomorphism,
\item a \dfn{fibration} if each $X(b)\ra Y(b)$ is a Kan fibration for
all $b\in \bool{B}$, and
\item a \dfn{weak equivalence} if $(L\Ex^\infty X)(b)\ra (L\Ex^\infty
Y)(b)$ is a weak equivalence of simplicial sets.  
\end{enumerate}
Here $\Ex^\infty\colon \Sh\bool{B}\ra \Psh\bool{B}$ denotes the
functor obtained by applying Kan's $\Ex^\infty$ functor at each
$b\in\bool{B}$.  This functor commutes with finite limits and
preserves fibrations; the same is true of the composite
$L\Ex^\infty\colon \Sh\bool{B}\ra\Sh\bool{B}$.

We say a map is a \dfn{trivial cofibration} if it is both a weak
equivalence and a cofibration, and we say a map is a \dfn{trivial
fibration} if it is both a weak equivalence and a fibration.  We say
an object $X$ is \dfn{fibrant} if the map $X\ra 1$ to the terminal
object is a fibration.

\begin{theorem}\label{thm-boolean-loc-model-cat}
With the above structure, $s\Sh\bool{B}$ is a closed model category in
the sense of Quillen.
\end{theorem}
\begin{proof}[Proof of (\ref{prop-local-fibs-are-global}) using
(\ref{thm-boolean-loc-model-cat})]
The cofibrations and weak equivalences in a closed model category
determine the fibrations.  Since $\Sh\bool{B}$ serves as its own
boolean localization, the weak equivalences (resp.\ fibrations) of
(\ref{thm-boolean-loc-model-cat}) are precisely the local weak
equivalences (resp.\ local fibrations) of
(\ref{subsec-model-cat-simplicial-sheaves}).  Thus the 
two model category structures coincide, and thus local fibrations are
model category theoretic fibrations.
\end{proof}

\subsection{Characterization of trivial fibrations}

\begin{lemma}
If $X\in s\Sh\bool{B}$ is fibrant, then $X(b)\ra (L\Ex^\infty X)(b)$ is
a weak equivalence for each $b\in\bool{B}$.
\end{lemma}
\begin{proof}
Since $X$ is fibrant, each map $X(b)\approx \prod_i X(b_i)\ra
\prod_i\Ex^\infty X(b_i)$ is a product of weak equivalences between
fibrant simplicial sets, and thus is a weak equivalence.  Thus the map
$$X(b)\ra (L\Ex^\infty X)(b)\approx
\coliml_{b=\amalg b_i}\prod_i \Ex^\infty X(b_i)$$
is a weak equivalence, since the colimit is taken over a directed set.
\end{proof}

\begin{corol}\label{cor-bool-loc-fib-we-crit}
If $X,Y\in s\Sh\bool{B}$ are fibrant, then $f\colon X\ra Y$ is a weak
equivalence if and only if $f(b)\colon X(b)\ra Y(b)$ is a weak
equivalence for each $b\in \bool{B}$
\end{corol}

The following lemma implies that if $f\colon X\ra Y\in s\Sh\bool{B}$
is a map such that each $f(b)$ is a weak equivalence of simplicial
sets, then $f$ is a weak equivalence.
\begin{lemma}\label{lemma-bool-loc-we-criterion}
If $f\colon X\ra Y\in s\Psh\bool{B}$ is a map of simplicial
presheaves such that each $f(b)\colon X(b)\ra Y(b)$ is a weak
equivalence, then $Lf\colon LX\ra LY$ is a weak equivalence of
simplicial sheaves.
\end{lemma}
\begin{proof}
First, note that if $W$ is a simplicial presheaf, then $L\Ex^\infty
LW\approx L\Ex^\infty W$.  Now if each $f(b)$ is
a weak equivalence, 
then so is each $\Ex^\infty f(b)$, and thus we conclude that
$L\Ex^\infty f$ is a weak equivalence using
(\ref{prop-sheafification-bool-loc}).  
\end{proof}

Let $y\colon \bool{B}\ra \Sh\bool{B}$ denote the canonical functor
sending $b$ to the representable sheaf $\hom_{s\Sh\bool{B}}(-,b)$.
Note that $y$ 
identifies $\bool{B}$ with the category of subobjects of the terminal
object $1$ in $\Sh\bool{B}$.

\begin{propos}\label{prop-bool-loc-tfib-char}
A map $f\colon X\ra Y\in s\Sh\bool{B}$ is a trivial fibration if and
only if each $f(b)\colon X(b)\ra Y(b)$ is a trivial fibration of
simplicial sets.
\end{propos}
\begin{proof}
The ``if'' part follows immediately from
(\ref{lemma-bool-loc-we-criterion}) and the definition of fibrations.

To prove the ``only if'' part, note that since $f$ is a fibration it
suffices to show that for each vertex $v\in Y_0(b)$ that the fiber of
$f(b)$ over $v$ is a contractible Kan complex.  Let $u\colon yb\ra Y$
be the map representing $v$, and form the pullback square
$$\xymatrix{
{P} \ar[r] \ar[d]_g 
& {X} \ar[d]^f
\\
{yb} \ar[r]^u
& {Y}
}$$
Note that $g$ is a fibration, and that the inclusion $yb\ra 1$ of
discrete objects is
easily seen to be a fibration.  Thus $P$ and $yb$ are fibrant.  Thus,
to show that $P(b)$ is contractible it suffices by
(\ref{cor-bool-loc-fib-we-crit}) to show that $g$ is a 
weak equivalence, since $yb(b)$ is a point.

The functor $L\Ex^\infty$ preserves pullbacks, fibrations, and weak
equivalences, and furthermore $L\Ex^\infty yb=yb$.  Thus 
$(L\Ex^\infty f)(b')$ is a trivial fibration for each $b'\in\bool{B}$,
whence so is each $(L\Ex^\infty g)(b')$, and thus $g$ is a weak
equivalence as desired.
\end{proof}

\subsection{Factorizations}

We produce factorizations of maps in $s\Sh\bool{B}$ by use of the ``small
object argument''.  

Choose an infinite cardinal $c>2^{|\bool{B}|}$ and let $\gamma$ be the
smallest 
ordinal of cardinality $c$.  
Then for each
$b\in\bool{B}$ the object 
$yb\in\Sh\bool{B}$ is \dfn{small} with respect to $\gamma$.  That is,
given a functor $X\colon \gamma\ra\Sh\bool{B}$, any map
$yb\ra\colim_{\alpha<\gamma} X_\alpha$
factors through some $X_\beta$ with $\beta<\gamma$.

\begin{lemma}\label{lemma-factor-cof-tfib}
Given $f\colon X\ra Y$, there exists a factorization $f=pi$ as a
cofibration $i$ followed by a trivial fibration $p$.
\end{lemma}
\begin{proof}
We inductively define a functor $X\colon \gamma\ra s\Sh\bool{B}$ as follows.  Let
$X(0)=X$.  Let $X(\alpha)=\colim_{\beta<\alpha}X(\beta)$ if
$\alpha<\gamma$ is a 
limit ordinal.  Otherwise, define $X(\alpha+1)$ by the pushout square
$$\xymatrix{
{\coprod yb\times\partial\Delta[n]} \ar[r] \ar[d]
& {\coprod yb\times\Delta[n]} \ar[d]
\\
{X(\alpha)} \ar[r]
& {X(\alpha+1)}
}$$
where the coproducts are taken over the set of all diagrams of the
form
$$\xymatrix{
{yb\times\partial\Delta[n]} \ar[r] \ar[d]
& {yb\times\Delta[n]} \ar[d]
\\
{X(\alpha)} \ar[r]
& {Y}
}$$
The desired factorization is $X\ra \colim_{\alpha<\gamma}X(\alpha)\ra Y$.
\end{proof}

\begin{lemma}\label{lemma-rlp-tfibs}
Trivial fibrations have the right lifting property with respect to
cofibrations. 
\end{lemma}
\begin{proof}
Note that
(\ref{prop-bool-loc-tfib-char}) and the ``choice'' axiom for
$\Sh\bool{B}$ (\ref{subsec-boolean-localization}) implies that trivial
fibrations are 
precisely the maps which have the right lifting property with respect
to all maps of the form $S\times \partial\Delta[n]\ra S\times
\Delta[n]$, where $S$ is any discrete object in $\Sh\bool{B}$.  Thus
the result follows when we note that if $i\colon 
A\ra B$ is a monomorphism in $s\Sh\bool{B}$, then we can write
$B_n\approx A_n\amalg S_n$ for some $S_n\in s\Sh\bool{B}$ since
$\Sh\bool{B}$ is boolean, and in this way construct an ascending
filtration $F_nB\subset B$ for $-1\leq n<\infty$ such that
$$F_nB\approx F_{n-1}B\bigcup_{S_n\times\partial\Delta[n]}
S_n\times\Delta[n],$$
$F_{-1}B\approx A$, and $\colim_n F_nB\approx B$.
\end{proof}

Let $\mathcal{C}$ denote the class of maps in $s\Sh\bool{B}$  which
are retracts of transfinite compositions of pushouts along maps of the
form $yb\times \Lambda^k[n]\ra yb\times\Delta[n]$, where
$b\in\bool{B}$ and $n\geq k\geq0$.

\begin{lemma}\label{lemma-factor-classc-fib}
Given $f\colon X\ra Y$, there exists a factorization $f=qj$ as a map
$j\in\mathcal{C}$ followed by a fibration $q$.
\end{lemma}
\begin{proof}
We perform the small object argument by taking pushouts along
coproducts of maps $yb\times 
\Lambda^k[n]\ra yb\times \Delta[n]$, indexed by diagrams of the form
$$\xymatrix{
{yb\times\Lambda^k[n]} \ar[r] \ar[d]
& {yb\times\Delta[n]} \ar[d]
\\
{X(\alpha)} \ar[r]
& {Y}
}$$
Otherwise, the proof is similar to that of (\ref{lemma-factor-cof-tfib}).
\end{proof}

\begin{lemma}\label{lemma-classc-is-tcof}
The class $\mathcal{C}$ is precisely the class of trivial cofibrations.
\end{lemma}
\begin{proof}
First we show that $\mathcal{C}$ consists of trivial cofibrations.  It
is already clear that every map in $\mathcal{C}$ is a cofibration.  

Suppose that $X$ is a simplicial sheaf.  Consider the following
pushout square in the category of simplicial \emph{presheaves}.
$$\xymatrix{
{\coprod yb\times\Lambda^k[n]} \ar[r] \ar[d]
& {\coprod yb\times \Delta[n]} \ar[d]
\\
{X} \ar[r]
& {Y}
}$$
That is, the products, coproducts, and pushouts are to be taken in the
category of presheaves.  It is then clear that $X(b')\ra Y(b')$ is a
weak equivalence for each $b'\in\bool{B}$, and thus by
(\ref{lemma-bool-loc-we-criterion}) the map $X\ra LY$ is a weak
equivalence in $s\Sh\bool{B}$; the sheaf $LY$ is the pushout of the
corresponding square of simplicial \emph{sheaves}.  Since the functor
$L\Ex^\infty$ 
commutes with directed colimits in $s\Sh\bool{B}$, and since weak
equivalences are 
closed under retracts, we may conclude that every map in $\mathcal{C}$
is a weak equivalence.

We show that any trivial cofibration $f\colon X\ra Y$ is in
$\mathcal{C}$ by a standard retract trick.  Namely, by
(\ref{lemma-factor-classc-fib}) we can factor $f=qj$ into a map
$j\in\mathcal{C}$ followed by a  fibration $q$.  But $q$ must also be a weak
equivalence since $f$ and $j$ are, and hence $q$ is a trivial
fibration.  Thus we can show using
(\ref{lemma-rlp-tfibs}) that $f$ is a retract of $j$ and thus
$f\in\mathcal{C}$. 
\end{proof}

\begin{proof}[Proof of (\ref{thm-boolean-loc-model-cat})]
Quillen's \cite{quillen-ratl-homotopy} axioms CM1, CM2, and CM3 are
clear: $s\Sh\bool{B}$ has all 
small limits and colimits, the classes of fibrations, cofibrations,
and weak equivalences are closed under retracts, and if any two of
$f$, $g$, and $gf$ are weak equivalences, then so is the third.
The factorization axiom CM5 follows from
(\ref{lemma-factor-cof-tfib}), (\ref{lemma-factor-classc-fib}), and
(\ref{lemma-classc-is-tcof}).  

One half of the lifting axiom CM4 follows from
(\ref{lemma-classc-is-tcof}), since fibrations are clearly
characterized by having the right lifting property with respect to all
maps of the form $yb\times \Lambda^k[n]\ra yb\times\Delta[n]$.  
The other half of CM4 is (\ref{lemma-rlp-tfibs}).
\end{proof}

\newcommand{\noopsort}[1]{} \newcommand{\printfirst}[2]{#1}
  \newcommand{\singleletter}[1]{#1} \newcommand{\switchargs}[2]{#2#1}

\end{document}